\documentstyle[11pt]{article}

\newcommand{\bbr}{\mbox{\boldmath $R$}}
\newcommand{\qed}{\hfill$\Box$}

\newcommand{\al}{\alpha}
\newcommand{\be}{\beta}

\newtheorem{thm}{Theorem}[section]
\newtheorem{pro}[thm]{Proposition}
\newtheorem{lem}[thm]{Lemma}
\newtheorem{cla}[thm]{Claim}

\newtheorem{obs}[thm]{Observation}

\newtheorem{exa}[thm]{Example}
\newtheorem{df}[thm]{Definition}
\newtheorem{rem}[thm]{Remark}

\def\sg{\sigma}

\begin{document}
\pagestyle{myheadings}

\title{Maximal degree in the Strong Bruhat Order of $B_n$}

\author{Tamar Seeman\thanks{Department of Mathematics, Bar-Ilan University, Ramat-Gan 52900}
\thanks{This research was partially supported by the Israel Science Foundation.}}

\date{September 11, 2006}

\maketitle{}

\bibliographystyle{acm}

\begin{abstract}

Given a permutation $\pi\in S_n$, let $\Gamma_-(\pi)$ be the graph on $n$ vertices 
$\{1,\ldots,n\}$
where two vertices $i<j$ are adjacent if $\pi^{-1}(i)>\pi^{-1}(j)$ and there are no integers 
$k$,
$i<k<j$, such that $\pi^{-1}(i)>\pi^{-1}(k)>\pi^{-1}(j)$. Let $\Gamma(\pi)$ be the graph 
obtained by
dropping the condition that $\pi^{-1}(i)>\pi^{-1}(j)$, i.e. two vertices are adjacent if the
rectangle $[i,\pi(i)]\times [j,\pi(j)]$ is empty. 
In the study of the strong order on permutation, Adin and Roichman introduced these graphs 
and computed their maximum number of edges.
We generalize these results to the Weyl group of signed permutations $B_n$, 
working with graphs on vertices $\{-n,\ldots,n\}\setminus \{0\}$, using new variants of 
a classical theorem of Tur\'an. 

\end{abstract}

\section{Introduction}

For a permutation $\pi$ in the symmetric group $S_n$, let the down degree be the number of
permutations covered by $\pi$ in the strong Bruhat order on $S_n$, and let the total degree 
be the valency of $\pi$ in the Hasse diagram of the strong Bruhat order.
Using a classical theorem of Tur\'an from
graph theory, Adin and Roichman showed that the down degree of 
a permutation in $S_n$ cannot exceed
$\lfloor n^2/4\rfloor$ \cite[Proposition 2.1]{AR}, and that the total degree cannot exceed
$\lfloor n^2/4\rfloor +n-2$ \cite[Theorem 3.1]{AR}. The authors also classified the 
extremal permutations, and computed the expected down degree of a random permutation. 

Geometrically, a permutation in $S_n$ can be plotted in $\bbr^2$ as a set of $n$ distinct points 
$\{(i,\pi(i))|~ 1\le i \le n\}$. In this context, the total degree of $\pi$ is equal to the
number of rectangles $[i,\pi(i)]\times [j,\pi(j)]$ containing no points in their interior. 
The
down degree is equal to the number of empty rectangles $[i,\pi(i)]\times [j,\pi(j)]$ with $i<j$ and $\pi(i)>\pi(j)$. 

In \cite{AFK}, Alon, { F$\ddot{{\rm u}}$redi and Katchalski extend 
the empty rectangle problem to $\bbr^d$, $d\ge 2$, giving asymptotic results for every fixed 
$d\ge 2$ as $n\rightarrow\infty$. Felsner~\cite{Fe} computes the number of empty rectangles in
$\bbr^3$, and gives different geometric applications. In a different direction, Sudakov \cite{Su} computes the maximal down and total degrees in a variant of the problem in which 
the number of points in the interior of a rectangle is bounded by some fixed integer $s$.  
For further results and applications see also~\cite{BB}. 

In this paper, we compute the maximal down degree and total degree in the group $B_n$ of 
signed permutations; that is, permutations which are anti-symmetric about the origin when 
plotted in $\bbr^2$. 
Section \ref{downBn} contains a description of the strong Bruhat order
on $B_n$ based on Coxeter generators (see also \cite{Inc}), and a corresponding definition of 
down degree. 

In section \ref{ddBn}, a variant
of Tur\'an's theorem is applied to show that the down degree of a permutation in $B_n$ cannot
exceed   
$\lfloor n^2/2 \rfloor$ (Proposition \ref{maxdownBn}), and a classification of the extremal 
permutations is given.  
Finally, the main result of section \ref{total} is Theorem \ref{maxBn}, which states that
the maximal total degree in $B_n$ is $\lfloor n^2/2 \rfloor+n-1$ for $n\ge 5$, and $4(n-1)$
for $n\le 5$. The proof of Theorem \ref{maxBn}, and a partial classification of the 
extremal permutations appears in section \ref{pftotal}.

\section{Down degree in $B_n$}\label{downBn}

For a positive integer $n$, let $[n]:=\{1,\ldots,n\}$.\\
Denote by $B_n$ the {\it hyperoctahedral group}, defined by
\[ B_n = \{\pi\in S_n:\ \pi(-i)=-\pi(i)\ i\in[n]\}.\]

\begin{df}
For a permutation $\pi\in B_n$ let the {\rm down degree} $d_-^B(\pi)$ be the number of
permutations in $B_n$ which are covered by $\pi$ in the strong Bruhat order. Let the 
{\rm up degree} $d_+^B(\pi)$ be the number of permutations which cover $\pi$ in this
order. The {\rm total degree} of $\pi$ is the sum
\[d^B(\pi)=d_-^B(\pi)+d_+^B(\pi).\]
\end{df}

Given integers $m$ and $n$, let $[m,n]:=\{m,\ldots,n\}$.\medskip 

For $a<b\in[n,n]\setminus 0$ let $t_{a,b}=t_{b,a}\in B_n$ be the transposition interchanging
$a$ and $b$, and let
\[ u_{a,b}=u_{b,a}=u_{\!-a,\!-b}=\left\{ {t_{a,b},~~~a=-b~~~~~~~~~~~~~~~~~~~~
}\atop{t_{a,b}t_{\!-a,\!-b},~~~ab>0~{\rm or}~
\pi^{-1}(a)\pi^{-1}(b)>0}\right.\]
If $a\ne -b$, $ab<0$ and $\pi^{-1}(a)\pi^{-1}(b)<0$, then $u_{a,b}$ is said to be undefined. 
\medskip

For $\pi\in B_n$ let
\[ \ell(\pi):={\rm min}\{k|\ \pi=s_{i_1}s_{i_2}\cdots s_{i_k}\} \]
be the {\it length} of $\pi$ with respect to the Coxeter generators 
$s_0=u_{1,-1}$, $s_i=u_{i,i+1},~i\in[1,n-1]$ of $B_n$.
Then
\begin{eqnarray*}
d_-^B(\pi)=\#\{u_{a,b}|~\ell(u_{a,b}\pi)=\ell(\pi)-1\}\\
d_+^B(\pi)=\#\{u_{a,b}|~\ell(u_{a,b}\pi)=\ell(\pi)+1\}\\
d^B(\pi)=d_-^B(\pi)+d_+^B(\pi)=\#\{u_{a,b}|~\ell(u_{a,b}\pi)=\ell(\pi)\pm 1\}
\end{eqnarray*}
We shall describe $\pi\in B_n$ by its sequence of values $[\pi(-n),\ldots,\pi(n)]$.\medskip

The following observation is based on \cite[Theorem 4.2]{Inc}.

\begin{obs}\label{coverB}
$\pi$ covers $\sg$ in the strong Bruhat order on $B_n$ if and only if there exist $i<k\in[-n,n]
\setminus 0$ such that
\begin{enumerate}
\item
$b:=\pi(i)>\pi(k)=:a$
\item
$u_{a,b}$ is defined
\item
$\sg=u_{a,b}\pi$
\item
There is no $i<j<k$ such that $a<\pi(j)<b$.
\end{enumerate}
\end{obs}

\begin{exa}
In $B_2$, $d[$-$2$-$112]=0$, $d[$-$1$-$221]=d[$-$21$-$12]=1$, and $d[$-$12$-$21]=d[1$-$22$-$1]
=d[2$-$11$-$2]=d[12$-$2$-$1]=d[21$-$1$-$2]=2$.
\end{exa}

\begin{df}
For $\pi\in B_n$ denote 
\[D_-^B(\pi):=\{u_{a,b}|\ \ell(u_{a,b}\pi)=\ell(\pi)-1\},\]
the strong $B$-descent set of $\pi$.
\end{df}

\begin{exa}
The strong $B$-descent set of $\pi=[$-$13$-$42$-$24$-$31]$ is
\[ D_-^B(\pi)=\{ u_{1,2}, u_{1,4}, u_{2,{\rm -}2}, u_{2,3}, u_{3,{\rm -}4}\}.\]
\end{exa}

\begin{pro}
The strong $B$-descent set $D_-^B(\pi)$ uniquely determines the permutation $\pi$.
\end{pro}

\noindent{\bf Proof.}
The proof is by induction on $n$, and is essentially the same as the proof of 
\cite[Proposition 1.8]{AR} in the $S_n$ case. 
Clearly the claim holds for $n=1$.

Let $\pi$ be a permutation in $B_n$, and let $\bar{\pi}\in B_{n-1}$ be the permutation
obtained by deleting $n$ and $-n$ from $\pi$. By the induction hypothesis $\bar{\pi}$ is
uniquely determined by the set 
\[ D_-^B(\bar{\pi})= D_-^B(\pi)\setminus\{u_{a,n}|~1\le a<n\}.\] Hence it suffices to 
determine the position of $n$ in $\pi$. 

Now, if $j:=\pi^{-1}(n)<n$ then $u_{\pi(j+1),n}\in D_-^B(\pi)$. Moreover, by Observation
\ref{coverB}, $u_{a,n}\in D_-^B(\pi)\Rightarrow a\ge\pi(j+1)$. Thus if 
$u_{-n,n}\in D_-^B(\pi)$ then $j=-1$. Similarly, if $u_{-n,n}\notin D_-^B(\pi)$, then 
$D_-^B(\pi)$ determines 
\[ \bar{\pi}(j)=\pi(j+1)={\rm min}\{a|\ u_{a,n}\in D_-^B(\pi)\} \] 
and therefore determines $j$. Note that this set of $a$'s is empty if and only 
if $j=n$. 
\qed\bigskip

\section{Maximal down degree in $B_n$}\label{ddBn}

In this section we prove

\begin{pro}\label{maxdownBn}
For every integer $n\ge 2$, 
\[ {\rm max}\{d^B_-(\pi)|\ \pi\in B_n\} = \lfloor n^2/2 \rfloor.\]
\end{pro}

\begin{df}
The {\rm underlying graph} of $\pi\in B_n$, denoted $\Gamma_-^B(\pi)$, is the graph whose set
of vertices is $[-n,n]\setminus 0$ and whose set of edges is
\[ \{\{a,b\}|\ u_{a,b}\in D_-^B(\pi)\}.\]
\end{df}

Define the function
\[ e:=e_\pi:\{(a,b)|\ a,b\in [-n,n]\setminus 0\}\longrightarrow \{0,1\}, \]
which maps pairs of adjacent vertices to $1$ and pairs of non-adjacent vertices to $0$. 
By definition, $e(a,a)=0$ for all $a$.\medskip

Denote $\pi[n]:= \{\pi(i)|~~i\in [n]\}.$

\begin{rem}\label{signs}
Let $a<b\in\pi[n]$.
If $e(a,-b)=1$, then $u_{a,-b}$ is defined so we must have $ab<0$.
\end{rem}

Let
\[ \al:=\al_\pi:\{(a,b)|\ a,b\in \pi[n]\}\longrightarrow \{0,1,2\}, \]
where \[ \al(a,b)= e(a,b) + e(a,-b).\]

Then
\begin{equation}\label{dpi1} 
d_-^B(\pi)=\sum_{a\le b\in\pi[n]}\al(a,b).
\end{equation}

\begin{df}\label{dega}
For vertex $a$ in $\Gamma_-^B(\pi)$, the degree of $a$, $d(a)$, is the number of
edges which are incident with $a$ in $\Gamma_-^B(\pi)$. 
Similarly, for vertices $a$ and $b$, $d(a\cup b)$ is the number of edges incident with 
either $a$ or $b$.
\end{df}

For the remainder of this section, let $a,b,c\in\pi[n]$.
\medskip

By Definition \ref{dega},
\begin{equation}\label{da}
d(a)= \sum_{b\in \pi[n]}\al(a,b), 
\end{equation}
and 
\begin{equation}\label{dab}
d(a\cup b)=d(a)+d(b)-\al(a,b).
\end{equation}

Given a set $S$ of vertices in $\pi[n]$, denote \[ -S=\{a~|-a\in S\}. \]

\begin{df}
For $S\subseteq\pi[n]$, $\Gamma_-^B(\pi)\setminus S$ is the subgraph 
of $\Gamma_-^B(\pi)$
consisting of all vertices not in $S$ or $-S$ and all edges not incident with a vertex in $S$
or $-S$. 
\end{df}
  
From Equation (\ref{dpi1}), we have 
\begin{equation}\label{dpi} 
d(\pi)=d(a)+d(\Gamma_-^B(\pi)\setminus a),~~~{\rm and}~~~
d(\pi)=d(a\cup b)+d(\Gamma_-^B(\pi)\setminus \{a,b\}).
\end{equation}

The proof of Proposition \ref{maxdownBn} requires the following four lemmas.

\begin{lem}\label{e}~\\
\noindent {\rm (a)}~~~$e(a,b)+e(a,c)+e(b,c)\le 2$,\\
\noindent {\rm (b)}~~~$e(a,-b)+e(a,-c)+e(b,-c)\le 2$.
\end{lem}

\noindent{\bf Proof.}
Let $a>b$, $e(a,b)=1$. If $e(a,c)=1$, then $c>a>b$ and $\pi^{-1}(c)<\pi^{-1}(a)<\pi^{-1}(b)$,
so $e(b,c)=0$. Conversely, $e(b,c)=1$ implies that $e(a,c)=0$, so (a) is proven. 

If $e(a,-b)=e(a,-c)=1$, then by Remark 
\ref{signs}, $ab<0$ and $ac<0$, hence $bc>0$ and $e(b,-c)=0$, implying (b).
\qed\bigskip

Lemma \ref{e} implies that 
\begin{equation}
\al(a,b)+\al(b,c)+\al(a,c)\le 4.
\end{equation}\label{le4} 
The following lemma, however, improves on this bound in certain cases.

\begin{lem}\label{a}
\noindent {\rm (a)}~~~If $\al(a,b)>0$, $\al(a,c)>0$ and $\al(b,c)>0$, then
\[\al(a,b)+\al(a,c)+\al(b,c)\le 3.\]
\noindent {\rm (b)}~~~If $a,b$ and $c$ are either all positive or
all negative, then \[\al(a,b)+\al(a,c)+\al(b,c)\le 2.\]  
\end{lem}

\noindent{\bf Proof.}
(a) Let $a>b$, $e(a,b)=1$. Then $\pi^{-1}(a)<\pi^{-1}(b)$. If $e(a,-c)=1$, then since
$\pi^{-1}(-c)<\pi^{-1}(a)$, we must have $-c>a>b$. But
$\pi^{-1}(c)<\pi^{-1}(a)<\pi^{-1}(b)$, so it follows that $e(b,-c)=0$. Conversely, if
$e(a,-c)=e(b,-c)=1$, then $e(a,b)=0$. Thus by Lemma \ref{e}(b), $e(a,-c)=e(b,-c)=1$ implies
that $\al(a,b)=0$. 


Now let $\al(x,y)>0$ for all $x<y\in\{a,b,c\}$. Then by the above argument,
there is at most one pair
$x<y\in\{a,b,c\}$ satisfying $e(x,-y)=1$. But by Lemma \ref{e}(a), at most two pairs
$x<y\in\{a,b,c\}$ satisfy $e(x,y)=1$, so we are done. 
\medskip

\noindent (b) Follows from Remark \ref{signs} and Lemma \ref{e}(a).
\qed\bigskip

\begin{lem}\label{a1}
\[\al(a,a)+\al(a,b)+\al(b,b)\le 2.\]
\end{lem}

\noindent{\bf Proof.}
Note that for every $x$, $e(x,x)=0$, so $\al(x,x)=e(x,-x)\le 1$. Also, if 
$x>0$ then $\al(x,x)=0$.

Without loss of generality, let $\pi^{-1}(a)<\pi^{-1}(b)$. 

Suppose that $\al(a,a)+\al(b,b)=2$. Then both $a<0$ and $b<0$, so by Remark 
\ref{signs}, $e(a,-b)=0$. Also, $\al(b,-b)=1$ implies that
$a<b$, which implies that $e(a,b)=0$. It follows that $\al(a,b)=0$. 

Now suppose that $\al(a,b)=2$. Then $e(a,b)=1$ implies that $a>b$, and
$e(a,-b)$ implies that $a<-b$, hence $b<0<a<-b$. But $\pi^{-1}(-b)<\pi^{-1}(a)<\pi^{-1}(b)$,
so it follows that $\al(b,b)=e(b,-b)=0$. Also, since $a>0$, $\al(a,a)=0$, so we are done.   
\qed\bigskip

\begin{lem}\label{uvwlem}
\[ \sum_{x\le y\in\{a,b,c\}}\al(x,y)\le 4.\]
\end{lem}\bigskip

\noindent{\bf Proof.}
For the lemma not to be satisfied, we would need either 
\[\al(a,b)+\al(a,c)+\al(b,c)=4,~~~{\rm or}~~~
\al(a,a)+\al(b,b)+\al(c,c)\ge 2.\] 

Suppose 
$\al(a,b)+\al(a,c)+\al(b,c)=4$. 
Then by Lemma \ref{a}(a), either $\al(a,b)=0$, $\al(a,c)=0$ or $\al(b,c)=0$.
For example, let $\al(a,b)=0$. Then $\al(a,c)=\al(b,c)=2$, so by Lemma
\ref{a1}, $\al(a,a)=\al(b,b)=\al(c,c)=0$ and we are done. 

Now suppose that $\al(a,a)+\al(b,b)+\al(c,c)\ge 2$. For example, let
$\al(a,a)=\al(b,b)=1$. Then by Lemma \ref{a1}, $\al(a,b)=0$, $\al(a,c)\le 1$ and 
$\al(b,c)\le 1$. If $\al(a,c)=\al(b,c)=1$, then by Lemma \ref{a1}, $\al(c,c)=0$.
Similarly, if $\al(c,c)=1$, then by Lemma \ref{a1}, $\al(a,c)=\al(b,c)=0$. Thus 
$\al(a,c)+\al(b,c)+\al(c,c)\le 2$, and we are done. 
\qed\bigskip

\noindent{\bf Proof of Proposition \ref{maxdownBn}.}
It is easy to see that the proposition is true for $n=2$. For $n=3$, the proof follows from
Equation (\ref{dpi1}) and Lemma \ref{uvwlem}. We prove the remaining cases by induction on $n$.
\medskip

Embed a signed permutation $\pi\in B_n$ as a permutation $\sg\in S_{2n}$ from 
$[-n,n]\setminus 0$ to itself. Let $G(\pi)$ be the graph on vertices $[-n,n]\setminus 0$, with 
edges between pairs of vertices $x,y$ such that $\sg$ covers $t_{x,y}\sg$ in the
strong Bruhat order on $S_{2n}$. By upper bound for $S_n$ (see \cite[Theorem 3.1]{AR}),
there are no more than $\lfloor (2n)^2/4\rfloor = n^2$ edges in this graph. 

Suppose $\al(a,a)=0$ for all $a\in\pi[n]$. 
Then every
covering relation $u_{x,y}$ in $\pi$ yields two edges in $G(\pi)$; namely, between $x$ and
$y$ and between $-x$ and $-y$. Thus $d^B_-(\pi)$ can be no more than 
$\lfloor n^2/2\rfloor$. 


Now suppose that $\al(a,a)=1$ for some $a$. Then by Lemma \ref{a1}, $\al(a,b)\le 1$ for all 
$b\ne a$, and therefore $d(a)\le n$. 
If $d(a)=n$, then for all $x\ne a$, $\al(a,x)=1$ and therefore $\al(x,x)=0$ by Lemma \ref{a1}. 
Let $b\ne a$, and suppose that $d(b)=n$. Since $\al(b,b)=0$, we must have $\al(b,c)=2$ for 
some $c$. But then $\al(a,c)+\al(a,b)+\al(c,c)=4$, contradicting Lemma \ref{a}(a).

It follows that $d(x)<n$ for some $x$. Thus by induction,
\[ d^B_-(\pi)=d(x)+d(\Gamma^B_-(\pi)\setminus x)< n + \lfloor (n-1)^2/2\rfloor
\le \lfloor n^2/2\rfloor. \]
\qed\bigskip

Next we classify the permutations with maximal down degree. For $\pi\in B_n$, we work only
with the vertices in $\pi[n]$.\bigskip

For $n=2$, the maximal permutations are $[$-$21]$, $[2$-$1]$, $[1$-$2]$, $[$-$2$-$1]$ and 
$[$-$1$-$2]$.

For $n=3$, the maximal permutations are $[$-$21$-$3]$, $[1$-$3$-$2]$, $[$-$31$-$2]$, 
$[2$-$3$-$1]$, $[$-$32$-$1]$,
$[3$-$2$-$1]$, $[3$-$21]$, $[12$-$3]$ and $[$-$2$-$1$-$3]$.

\begin{lem}\label{incr}
Let $n\ge 4$, $\pi\in B_n$ a permutation with maximal down degree. 
Then $\pi$ consists of a shuffle of an increasing 
sequence of positive integers and an increasing sequence of negative integers.
\end{lem}

\noindent{\bf Proof.}
Suppose $\al(a,b)=1$ for some $a\ne b$. As shown in the proof of Proposition \ref
{maxdownBn}, 
\[ d^B_-(\pi)= d(a\cup b)+d(\Gamma^B_-(\pi)\setminus\{a,b\}),\]
where
\[ d(a\cup b)=\al(a,a)+\al(a,b)+\al(b,b)+\sum_{c\notin\{a,b\}}(\al(a,c)+\al(b,c))\le 2n-2,\]
and \[ d(\Gamma^B_-(\pi)\setminus\{a,b\})\le \lfloor(n-2)^2/2\rfloor.\]
Thus $\pi$ maximal implies that both \[ d(a\cup b)=2n-2~~~{\rm and}~~~ 
d(\Gamma^B_-(\pi)\setminus\{a,b\})= \lfloor(n-2)^2/2\rfloor.\]

By Lemma \ref{e}(a),
\[\sum_{c\notin\{a,b\}}(\al(a,c)+\al(b,c))\le 2(n-2),\]
hence $d(a\cup b)=2n-2$ requires that $\al(a,a)+\al(b,b)=1$. Since $\al(x,x)=0$ for all $x>0$,
it follows that either $a<0$ or $b<0$.
Thus $\al(x,y)=0$ for all pairs of positive integers $x$ and $y$, so the positive
integers must be increasing.

Now suppose that $a<0$ and $b<0$. 
If all of the other $n-2$ integers are positive, then
by Lemma \ref{a}(b), the subgraph $\Gamma^B_-(\pi)\setminus\{a,b\}$ is triangle-free, hence by
Tur\' an's theorem 
\[ d(\Gamma^B_-(\pi)\setminus\{a,b\})\le\lfloor (n-2)^2/4\rfloor<\lfloor (n-2)^2/2\rfloor,\]
so $\pi$ is not maximal.    
We therefore assume that $c<0$ for some $c\notin\{a,b\}$. But then by Lemma \ref{a}(b), 
$\al(a,c)+\al(b,c)\le 1$, implying that $d(a\cup b)<2n-2$. It follows that $\al(x,y)=0$ 
for all pairs of negative integers $x$ and $y$, so the negative integers must be
increasing.
\qed\bigskip

\begin{pro}
Let $n\ge 4$, $\pi\in B_n$, $d_-^B(\pi)=\lfloor n^2/2\rfloor$. Then 
\[ \pi=[1,2,\ldots,m,-n,-(n-1),\ldots,-(m+1)], \]
where $m\in\{\lfloor n/2\rfloor, \lceil n/2\rceil\}.$
\end{pro}

\noindent{\bf Proof.}
Suppose that $\al(a,a)=1$ for some $a$. Then $a<0$, and by Lemma \ref{a1}, $\al
(a,b)\le 1$ for all $b\ne a$, so by Equation (\ref{da}), $d(a)\le n$.  

If $a$ is the only negative integer, then by
Lemma \ref{incr}, $\pi\setminus a$ consists of an increasing sequence of positive integers,
hence $d(\Gamma^B_-(\pi)\setminus a)=0$, implying that $d(\pi)\le n$ and $\pi$ is not 
maximal. On the other hand, if there are two negative integers $b$ and $c$ other than $a$, then
by Lemma \ref{incr}, $\al(a,b)=\al(a,c)=0$, which implies that $d(a)\le n-2$. But then 
\[ d^B_-(\pi)=d(a)+d(\Gamma^B_-(\pi\setminus a))\le n-2 +\lfloor(n-1)^2/2\rfloor
<\lfloor n^2/2\rfloor,\] and again $\pi$ is not maximal. 

Suppose, however, that there is exactly one other negative integer $b\ne a$. 
Then $d(\Gamma^B_-(\pi\setminus \{a,b\}))=0$, so by Equation (\ref{dpi}),
$d_-^B(\pi)=d(a\cup b).$
But $d(a\cup b)=d(a)+d(b)\le n-2+2(n-2)=3n-6<\lfloor n^2/2\rfloor,$ implying that $\pi$ is not
maximal.

It follows that $\pi$ maximal implies that $\al(a,a)=0$ for all $a$, and therefore by 
Lemma \ref{incr}, 
\[ d^B_-(\pi)= \sum_{a>0,~b<0}\al(a,b)\le 2k(n-k), \]
where $k$ is the number of positive integers in $\pi$. Thus for $\pi$ to be maximal, we must
have $k\in\{\lfloor n/2\rfloor, \lceil n/2\rceil\},$ and $\al(a,b)=2$ for all $a>0$, $b<0$.
But for all $a>0$ and $b<0$, $\al(a,b)=e(a,b)+e(a,-b)=2$ implies that 
both $\pi^{-1}(a)<\pi^{-1}(b)$ and $a<-b$, so the proof follows.
\qed\bigskip

\section{Total Degree in $B_n$}\label{total}

Proposition \ref{maxdownBn} holds for $d_+^B$ instead of $d_-^B$. Thus the maximal value of 
the total degree $d^B=d^B_-+d^B_+$ cannot exceed $2\lfloor n^2/2\rfloor=
\lfloor n^2\rfloor$. For $n\ge 4$, the maximal value is smaller.


\begin{thm}\label{maxBn} 
\[ {\rm max}\{d^B_-(\pi)|\ \pi\in B_n\} = \left\{{4(n-1),~~2\le n\le 5,}\atop 
{\lfloor n^2/2 \rfloor+n-1,~~n\ge 5}\right. \]
\end{thm}

\begin{df}
The {\rm underlying graph} of $\pi\in B_n$, denoted $\Gamma^B(\pi)$, is the graph whose set
of vertices is $[-n,n]\setminus 0$ and whose set of edges is
\[ \{\{a,b\}|\ \ell(u_{a,b}\pi)-\ell(\pi)=\pm 1\}.\]
\end{df}

Define the edge function $e_\pi$ on $\Gamma^B(\pi)$,
which maps pairs of adjacent vertices to $1$ and pairs of non-adjacent vertices to $0$. 

\begin{rem}\label{tsigns}
For $a,b\in\pi[n]$, $e(a,-b)=1$ implies that $ab<0$.
\end{rem}

Let $\be_\pi:\{(a,b)|\ a,b\in [n]\}\longrightarrow \{0,1,2\},$
where \[ \be_\pi(a,b)= e_\pi(a,b) + e_\pi(a,-b).\]
Then
\begin{equation}\label{tdpi1} 
d^B(\pi)=\sum_{a\le b\in\pi[n]}\be_\pi(a,b).
\end{equation}





Similar to the case of down degree, we have
\begin{equation}\label{tda}
d(a)= \sum_{b\in\pi[n]}\be_\pi(a,b), 
\end{equation}
and 
\begin{equation}\label{tdab}
d(a\cup b)=d(a)+d(b)-\be_\pi(a,b).
\end{equation}

\begin{df}
For $S\subseteq\pi[n]$, $\Gamma^B(\pi)\setminus S$ is the subgraph 
of $\Gamma^B(\pi)$
consisting of all vertices not in $S$ or $-S$ and all edges not incident with a vertex in $S$
or $-S$. 
\end{df}
  
Equation (\ref{tdpi1}) implies that for $a,b\in\pi[n]$,
\begin{equation}\label{tdpi} 
d(\pi)=d(a)+d(\Gamma^B(\pi)\setminus a),~~~{\rm and}~~~
d(\pi)=d(a\cup b)+d(\Gamma^B(\pi)\setminus \{a,b\}).
\end{equation}

\begin{df}\label{piS}
Let $\pi\in B_n$, $S\subseteq \pi[n]$, $|S|$ the size of $S$. Then $\pi\setminus S$ is the permutation in $B_{n-|S|}$ obtained by eliminating from $\pi$ every integer in $S$ or $-S$.
\end{df}

\begin{exa}
Let $\pi=3$-$1$-$2$, then $d(\pi\setminus$-$1 )=d(3$-$2)=4$. On the other hand,
$d(\Gamma^B(\pi)\setminus$-$1)=3$.
\end{exa}

Note that in the above example, $d(\pi\setminus$-$1)> d(\Gamma^B(\pi)\setminus$-$1)$.

\begin{df}
For $S\subseteq\pi[n]$, 
\begin{equation}\label{rS}
r(S)=d(\pi\setminus S)-d(\Gamma^B(\pi)\setminus S)=
\sum_{a,b\in\pi[n]\setminus S}(\be_{\pi\setminus S}(a,b)-\be_\pi(a,b))
\end{equation}
By definition, $r(S)\ge 0$ for every set $S$.
\end{df}


Equations (\ref{tdpi}) and (\ref{rS}) together give 
\begin{equation}\label{da-ra}
d(\pi)=d(\pi\setminus a)+d(a)-r(a),
\end{equation}
 and
\begin{equation}\label{dab-rab}
d(\pi)=d(\pi\setminus a)+d(a\cup b)-r(\{a,b\}).
\end{equation}
\medskip

A permutation $\pi\in B_n$ can be represented geometrically by plotting the points 
$(i,\pi(i))$, $i\in[-n,n]\setminus 0$ in the Cartesian plane. The rectangle with 
$(i,\pi(i))$ and $(j,\pi(j))$ at opposite corners is denoted  
\[ [i\times j] := [i,\pi(i)]\times[j,\pi(j)]. \]
By Observation \ref{coverB}, two vertices
$\pi(i)$ and $\pi(j)$ are adjacent in $\Gamma^B(\pi)$ if and only if 

\noindent (a)~~~the rectangle $[i\times j]$ is empty, and

\noindent (b)~~~either $j=-i$ or $(0,0)\notin [i\times j]$.\medskip

\noindent Thus by Equation (\ref{tdpi1}),  
$d(\pi)$ is the number of empty rectangles $[i\times j]$ in $\pi$, such that $i\in[n]$, 
$|j|\ge i$, and conditions (a) and (b) are satisfied.
Note that $\pi(i)$ and $\pi(j)$ are adjacent if and only if $\pi(-i)$ and $\pi(-j)$
are adjacent, so we need not count both $[i\times j]$ and $[-i\times -j]$. Hence the requirement that $i\in[n]$ and $|j|\ge i$. 

Given a set $S$ of vertices in $\pi[n]$, $r(S)$ is the number of {\it nonempty} rectangles 
$[i\times j]$ satisfying (b), and with the following property: For every $a\in [i\times j]$, either $a\in S$ or $-a\in S$. 

\section{Proof of Theorem \ref{maxBn}}\label{pftotal}


\begin{lem}\label{maxBnlem}
Let $n\ge 4$. Then either\\ 
\noindent {\rm (a)}~~$d(a)-r(a)\le n$ for some $a\in \pi[n]$, or\\
\noindent {\rm (b)}~~$d(a\cup b)-r(\{a,b\})\le 2n$ for some $a,b\in \pi[n]$.
\end{lem}

\noindent{\bf Proof.}
Let $-\pi$ be the permutation obtained from $\pi$ by multiplying each integer by $-1$.
Then $d(-\pi)=\pi$, so without loss of generality we can assume that $1\in\pi[n]$.\medskip

{\it Case 1: $0<\pi^{-1}(2)<\pi^{-1}(1)$.}

Let $a\in\pi[n]$, $\be_\pi(1,a)=2$. Then by Remark \ref{tsigns}, $a<-2$. But then $2\in[1\times -a]$, a contradiction. It follows that $d(1)\le n$. 
\medskip

{\it Case 2: $0<\pi^{-1}(-2)<\pi^{-1}(1)$.}

Let $a\in\pi[n]$, $\be_\pi(1,a)=2$. Then by Remark \ref{tsigns}, 
$a\le -2$. If $a$ is left of $-2$, then 
$-2\in [a\times 1]$, a contradiction. Similarly, if $a$ is right of 
$-2$, then $2\in [-a\times 1]$, a contradiction. It follows that
$-2$ is the only possible value for $a$, hence $d(1)\le n+1$.

If $\be_\pi(1,x)=0$ for some $x$, then $d(1)\le n$ and we are done. Therefore assume that 
$\be_\pi(1,x)\ge 1$ for all $x$, so the positive integers left of $1$ are increasing.

If $d(-2)\le n+1$, then by Equation (\ref{tdab}),
$d(1\cup -2)\le 2n$, and we are done. Therefore
assume that $d(-2)\ge n+2$, so for at least one integer $a\ne 1$, $\be_\pi(-2,a)=2$. Choose $a$ to be the
rightmost such integer. By Remark \ref{tsigns}, $a>1$, implying that $a$ appears left of $1$. Suppose 
that a positive integer $x$ appears left of $a$. Since the positive integers left of $1$ are
increasing, $x\in[2\times a]$, contradicting the assumption that $\be_\pi(-2,a)=2$. It
follows that $a$ is the lowest and leftmost positive 
integer left of $1$, and
for every $w\in \pi[n]\setminus \{1,a\}$, $\be_\pi(-2,w)\ne -2$. Thus
$d(-2)\le n+2$, and by Equation (\ref{tdab}), $d(1\cup -2)=(n+1)+(n+2)-2=2n+1$.

Thus the lemma is proven if at least one of the following conditions is satisfied:

\noindent (A) $r(1) \ge 1$

\noindent (B) $r(-2) \ge 2$

\noindent (C) $d(-2)\le n+1$ (which implies $d(1\cup-2)\le 2n$)

\noindent (D) $r(\{1,-2\})\ge 1$.

Suppose there is a positive integer right of $1$. Let $u$ be the leftmost such integer. Then
$1\in [-2\times u]$, and condition (A) is satisfied. Therefore assume $1$ is the rightmost 
positive integer. Now suppose there is a negative integer $v$ right of $1$.
Then for all positive $z$ left of $1$, $1\in [z\times v]$, hence $\be_\pi(z,v)\le 1$. Also, since $2\in [1\times -v]$, $\be_\pi(1,v)\le 1$. 
Thus $d(v)\le n$ and we are done. Therefore assume that $1=\pi(n)$.\medskip

{\it Case 2a: $\pi^{-1}(-3)>0$.}
 
We assume that $\be_\pi(-2,a)=2$ for some $a>2$, since otherwise condition (C) would be
satisfied. Thus by assumption $3\notin [2\times a]$, so $-3$ must 
appear to the right of $-2$. It follows that $-2$ and $2$ are the only integers inside
$[3\times -3]$, hence $r(-2)\ge 1$.

Recall that the integer $a>2$ with $\be_\pi(-2,a)=2$ must be the lowest and leftmost positive 
integer left of $1$.
If $a$ is left of $-2$, then $-2\in [a\times -3]$, and condition
(B) is satisfied. Suppose, however, that $a$ is right of $-2$. If there is another positive 
integer $x\ne a$ appearing left of $1$, then $x>a$ and $x$ appears right of $a$. But then
$a$ appears in both $[-2\times x]$ and 
$[2\times x]$, hence $\be_\pi(-2,x)=0$, so condition (C) is satisfied. It follows that if $a$ is right of $-2$, then $a$ is the only positive integer other than $1$. 

Now $2\in [1\times 3]$, so $\be_\pi(1,-3)\le 1$. Since $a$ is the only positive integer other 
than $1$, $\be_\pi(-3,w)\le 1$ for all $w\in\pi[n]\setminus a$.  
Also, recall that $-2\in [-3\times 3]$, so $\be_\pi(3,3)=0$.
It follows that $d(-3)\le n$ and we are done.\medskip

{\it Case 2b: $\pi^{-1}(3)>0$.}

Recall that by assumption $1=\pi(n)$ and $\be_\pi(1,w)\ge 1$ for all $w\in\pi[n]$. Suppose some
$x<-3$ appears left of $-2$. 
Then $-2\in [1\times x]$ and $3\in [1\times -x]$, so $\be_\pi(1,x)=0$, a contradiction.  
Therefore by assumption there are no negative integers left of $-2$, so any integers 
appearing left of $-2$ must be positive (and in increasing order).\medskip

(i) Suppose $3$ is left of $-2$. Then $3=\pi(1)$. 
Let $x<-3$, so by assumption, $x$ appears right of $-2$. Then $-2\in [3\times x]$, implying 
that $\be_\pi(3,x)\le 1$. It follows that $\be_\pi(3,w)\le 1$ for all $w\in\pi[n]\setminus -2$,
and therefore $d(3)\le n+1$.

If $|\pi(2)|>3$, then $3$ and $-3$ are the only integers in $[\pi(-2)\times \pi(2)]$, hence 
$d(3)-r(3)\le n$ and we are done.
Suppose, however, that $\pi(2)=-2$. Recall that $d(-2)\le n+2$. Since $\be_\pi(3,-2)=2$, we therefore have $d(3\cup -2)\le (n+1) + (n+2)-2=2n+1$. Thus the lemma
is satisfied if we can show that $r(\{3,-2\})\ge 1$. But since $n\ge 4$ and $1=\pi(n)$, it 
follows that $|\pi(3)|>3$, and $3,-2,-3$ and $2$ are the only integers inside 
$[\pi(-3)\times \pi(3)]$, so we are done.\medskip 

(ii) Suppose $3$ is right of $-2$. If an additional positive integer $x>3$ appears left of $1$,
then since the positive integers are in increasing order, $x$ must be right of $3$.
But then $3\in [-2\times x]$ and $3\in [2\times x]$, so $\be_\pi(2,x)=0$ and
condition (C) is satisfied. Therefore assume that $1$ and $3$ are the only positive integers. 
Since $n\ge 4$, it follows that there appears some integer $z<-3$. 

Recall that there are no negative integers left of $-2$. Thus $z$ must be right of $-2$,
so $2\in [-z\times z]$. Also,
recall that $\be_\pi(1,w)\le 1$ for all $w\ne -2$. Thus $\be_\pi(1,z)\le 1$,
so for all $w\ne 3$, $\be_\pi(w,z)\le 1$. It follows that $d(z)\le n$ and we are done. 
\medskip

{\it Case 3: $2$ or $-2$ appears right of $1$.}

The integers $w$ satisfying $\be_\pi(1,w)=2$ are in increasing order, with no more than one 
appearing to the left of $1$. 
If $r(1)$ is greater than or equal to the number of such integers, then 
$d(1)-r(1)\le n$ and we are done. We assume, however, that this is not the case, and therefore 
one can find an integer $w$ satisfying $\be_\pi(1,w)=2$ and the
following property: For every integer $x$, if $1\in [w\times x]$, then some $y$ such that 
$|y|\ne 1$ is also in $[w\times x]$.
\medskip
 
{\it Case 3a: $w$ is right of $1$.}

We show that   
\begin{equation}\label{1w}
d(1\cup w)-r(\{1,w\})\le 2n.
\end{equation}

Now $1$ and $-1$ are inside $[w\times -w]$, so by assumption
the rectangle $[w\times -w]$ contains in its 
interior some integer $y$ other than $1$ and $-1$. In other words, there appears some $y$ left 
of $w$ such that $1<|y|<-w$. Now, if $y<0$ then $-y\in [1\times -w]$, contradicting the
assumption that $\be_\pi(1,w)=2$. Also, if $y>0$ and $y$ appears left of $1$, then 
$y\in [1\times -w]$, again a contradiction. It follows that $y$ is positive and appears between 
$1$ and $w$. Choose $y$ to be the leftmost integer right of $1$ such that $1<y<-w$. 


Consider the integers $x$ satisfying 
\begin{equation}\label{be3}
\be_\pi(w,x)+\be_\pi(1,x)=3.
\end{equation} 
Let $x>0$, so for $x$ to satisfy (\ref{be3}), we need $\be_\pi(w,x)=2$ and $\be_\pi(1,x)=1$.
Suppose $x$ appears right of $w$. If $x>y$, then $y\in [1\times x]$, so 
$\be_\pi(1,x)=0$. On the other hand, if $x<y$, then $y\in [x\times -w]$, and therefore 
$\be_\pi(w,x)\le 1$. It follows that to satisfy Equation (\ref{be3}), $x$ must appear left of
$w$. More precisely, $x$ must appear between $1$ and $w$, since $x$ left of $1$ would imply
that $1\in[w\times x]$.

Now suppose that a second positive integer $z$ satisfies Equation (\ref{be3}), and assume 
that $z$ appears between $x$ and $w$. If $z<x$, then $z\in[w\times x]$ and therefore
$\be_\pi(w,x)\le 1$, a contradiction. Similarly, if $z>x$, then $x\in[1\times z]$, and
therefore $\be_\pi(1,z)=0$, again a contradiction. It follows that at most one positive 
integer satisfies Equation (\ref{be3}).

Now let $x<0$, so for $x$ to satisfy (\ref{be3}), we need $\be_\pi(1,x)=2$ and 
$\be_\pi(w,x)=1$. Recall that the integers $x$ satisfying $\be_\pi(1,x)=2$ form an increasing 
sequence which includes $w$. Thus at most two integers in this sequence also satisfy 
$\be_\pi(w,x)=1$. It follows that at most two negative integers satisfy Equation (\ref{be3}).
Altogether, no more than three integers satisfy Equation (\ref{be3}).

Since $\be_\pi(1,1)=1$ but $\be_\pi(w,w)=0$,  
it follows that
$d(1)+d(w)\le 2n+4$, which implies that $d(1\cup w)\le 2n+2$. Note that $r(\{1,w\})\ge 1$, since $1$ and $-1$ are the only integers inside $[y\times -y]$. 
(Follows from the fact that $\be_\pi(1,w)=2$ and $y<-w$.)

Now if only one negative integer satisfies Equation (\ref{be3}), then $d(1\cup w)\le 2n+1$, 
which implies (\ref{1w}), and we are done. On the other hand, if two negative integers 
$u$ and $v$ satisfy Equation (\ref{be3}), then $w$ appears between $u$ and $v$ in the 
increasing sequence of integers satisfying $\be_\pi(1,x)=2$, and therefore $w$ is the only nteger in $[u\times v]$, hence $r(\{1,w\})\ge 2$, again implying (\ref{1w}) so we are done. 
\medskip

{\it Case 3b: $w$ is left of $1$}

Let $s$ be the number of integers $x$ satisfying $\be_\pi(1,x)=2$. Recall that these $s$ 
integers are in increasing order, with no more than one (namely $w$) appearing to the left of
$1$. 
Assume that for every $x$ right of $1$ with $\be_\pi(1,x)=2$, $1$ and $-1$ 
are the only integers inside $[x\times -x]$, since otherwise we could apply Case 3a. 
Thus $r(1)\ge s-1$. 

If $2$ appears right of $1$,
then $1$ is the only integer inside $[w\times 2]$, so $r(1)\ge s$ and we are done. Otherwise,
$-2$ appears right of $1$. If there is a positive integer left of $1$, then denoting by $y$ the 
rightmost such positive integer, $1$ is the only integer in $[y\times -2]$, hence again
$r(1)\ge s$ and we are done. Suppose, however,
that no positive integers appear left of $1$. Then $\be_\pi(w,x)\le 1$ for all $x\in\pi[n]\setminus 1$, and therefore $d(w)\le n+1$. We assume that $\be_\pi(w,x)>0$ for all $x$, since otherwise we would have $d(w)\le n$.    

Since $w$ is the leftmost integer in the increasing sequence of $s$ integers satisfying
$\be_\pi(1,x)=2$, each of the $s-2$ rightmost integers satisfies $\be_\pi(w,x)=0$. But by assumption $\be_\pi(w,x)>0$ for all $x$, so it follows that $s=2$ and therefore 
$\be_\pi(1,x)\le 1$ for all $x\in\pi[n]\setminus\{-2,w\}$. Thus $d(1)\le n+2$, which implies that  
$d(1\cup w)\le 2n+1$. But since $-2$ appears to the right of $1$, $1$ and $-1$ are the only 
integers inside $[2\times -2]$, hence $r(\{1,w\})\ge 1$ and we are done.
\qed\bigskip

\noindent{\bf Proof of Theorem \ref{maxBn}}
Proposition \ref{maxdownBn} holds for $d_+^B$ instead of $d_-^B$. Thus the maximal value of 
the total degree $d^B=d^B_-+d^B_+$ cannot exceed $2\lfloor n^2/2\rfloor=
\lfloor n^2\rfloor$. The proof follows for $n=2,3$. The case $n=5$ was verified by computer. 
We prove the remaining cases by induction. 

Let $n\ge 4$, $n\ne 5$. Suppose that $d(a)-r(a)\le n$ for some $a\in \pi[n]$.  
Then by Equation (\ref{da-ra}), $d(\pi)\le d(\pi\setminus a)+n$, where 
$\pi\setminus a\in B_{n-1}$. If $n=4$, then the proof follows easily by induction. Similarly,
if $n\ge 6$, then  
\[ d(\pi)\le d(\pi\setminus a)+n \le \lfloor (n-1)^2/2\rfloor +n-2 + n\le
\lfloor n^2/2\rfloor +n-1, \]
so again the proof follows by induction. 

On the other hand, if $d(a)-r(a)>n$ for all $a\in\pi[n]$, then by Lemma \ref{maxBnlem},
$d(a\cup b)-r(\{a,b\})\le 2n$ for some $a,b\in \pi[n]$.
Thus by Equation (\ref{dab-rab}), $d(\pi)\le d(\pi\setminus \{a,b\})+2n$, where 
$\pi\setminus \{a,b\}\in B_{n-2}$. If $n=4$ or $6$, then the proof follows by induction. Similarly, 
if $n> 6$, then  
\[ d(\pi)\le d(\pi\setminus \{a,b\})+2n \le \lfloor (n-2)^2/2\rfloor +n-3 + 2n\le
\lfloor n^2/2\rfloor +n-1, \]
which implies the proof. 
\qed\bigskip

Next we classify the permutations with maximal total degree. For $\pi\in B_n$, we work only
with the vertices in $\pi[n]$. 

\begin{cla}
Let $\pi$ be a permutation of maximal total degree.
\begin{enumerate}
\item If $n=2$, $\pi\in\{[2$-$1]$, $[$-$21]\}$.
\item If $n=3$, $\pi\in\{[$-$32$-$1]$, $[3$-$21]\}$.
\item If $n=4$, $\pi\in\{[4$-$32$-$1]$, $[$-$43$-$21]$, $[4$-$3$-$21]$, $[$-$432$-$1]\}$.
\item If $n=5$, there are $112$ maximal permutations.
\item If $n\ge 6$, then the number of maximal permutations is $8$ for even $n$ and $16$ 
for odd $n$. The extremal permutations have one of the following forms:
\[ \pi_0:=[1,2,\ldots,m,{\rm -}n,{\rm -}(n-1),\ldots{\rm -}(m+1)]~~~~(m\in\{\lfloor 
n/2\rfloor,\lceil n/2\rceil\}),\]
and the permutations obtained from $\pi_0$ by one or more of the following operations:
\[ \pi \mapsto -\pi~~~~({\rm multiplying~every~integer~by~-}1),~~~~~~~~~~~~~~~~~~~~~~~~~~~\]
\[ \pi \mapsto u_{m,m+1}\pi
~~~~({\rm interchanging~} m~{\rm and~}m+1,~{\rm if~}u_{m,m+1}~{\rm is~defined}),\]
\[ \pi \mapsto u_{m,{\rm -}n}\pi~~~~({\rm interchanging~} m~{\rm and~-}n,~{\rm if~}
u_{m,{\rm -}n}~{\rm is~defined}).~~~~~~~\]
\end{enumerate}
\end{cla}

\noindent{\bf Proof.}
By computer verification.
\qed

\end{document}